\documentclass[12pt]{amsart}
\usepackage{amscd,amssymb,graphics,color,a4wide}

\usepackage{mathrsfs}

\newtheorem{theorem}{Theorem}[section]
\newtheorem{lemma}[theorem]{Lemma}
\newtheorem{proposition}[theorem]{Proposition}

\theoremstyle{definition}
\newtheorem{definition}[theorem]{Definition}
\newtheorem{example}[theorem]{Example}
\newtheorem{examples}[theorem]{Examples}

\theoremstyle{remark}
\newtheorem{remark}[theorem]{Remark}

\newcommand{\NN} {\mathbb{N}}
\newcommand{\ZZ} {\mathbb{Z}}

\newcommand{\RR} {\mathbb{R}}
\newcommand{\CC} {\mathbb{C}}

\newcommand{\PP} {\mathbb{P}}

\newcommand {\shExt} {\mathcal{E} \!\text{\textit{xt}}}

\newcommand {\shL}  {\mathcal{L}}
\newcommand {\shM}  {\mathcal{M}}

\newcommand {\shN}  {\mathcal{N}}

\newcommand {\shT}  {\mathcal{T}}

\newcommand {\shP}  {\mathcal{P}}
\newcommand {\shU}  {\mathcal{U}}

\newcommand {\shX}  {\mathcal{X}}

\newcommand {\Aff}  {\operatorname{Aff}}

\newcommand {\dual} {\vee}

\newcommand {\Hom}  {\operatorname{Hom}}

\newcommand {\im}  {\operatorname{im}}

\newcommand {\M} {\mathscr{M}}

\renewcommand{\O}  {\mathcal{O}}

\renewcommand{\P}  {\mathscr{P}}

\newcommand {\Spec} {\operatorname{Spec}}

\newcommand {\T} {\shT}
\newcommand {\X} {\shX}

\newcommand {\U} {\mathscr{U}}

\def\mydate{\ifcase\month \or January\or February\or March\or
April\or May\or June\or July\or August\or September\or October\or 
November\or December\fi \space\number\day,\space\number\year}

\begin{document}
\def\mapright#1{\smash{
  \mathop{\longrightarrow}\limits^{#1}}}
\def\exact#1#2#3{0\rightarrow#1\rightarrow#2\rightarrow#3\rightarrow0}
\def\mapup#1{\Big\uparrow
   \rlap{$\vcenter{\hbox{$\scriptstyle#1$}}$}}
\def\mapdown#1{\Big\downarrow
   \rlap{$\vcenter{\hbox{$\scriptstyle#1$}}$}}
\def\dual#1{{#1}^{\scriptscriptstyle \vee}}
\input epsf.tex

\title{Affine Manifolds,
Log Structures, and Mirror Symmetry}

\author[GROSS and SIEBERT]{Mark Gross, Bernd Siebert} 
\address{Mathematics Institute,
University of Warwick, Coventry CV4 7AL, United Kingdom}
\address{Department of Mathematics,
UCSD, La Jolla, CA 92093-0112, United States}
\email{mgross@maths.warwick.ac.uk}
\thanks{This work was partially supported by NSF grant 0204326.}

\address{Mathematisches Institut,
 Albert-Ludwigs-Universit\"at Freiburg, Eckerstrasse 1,
D-79104 Freiburg, Germany}
\email{bernd.siebert@math.uni-freiburg.de}

\maketitle
\bigskip
\section*{Introduction.}

Mirror symmetry between Calabi-Yau manifolds is inherently about
degenerations: a family $f:\X\rightarrow S$ of Calabi-Yau varieties
where $S$ is a disk, $\X_t$ is a non-singular Calabi-Yau
manifold for $t\not=0$, and $\X_0$ a singular variety. Much information
about the singular fibre is carried in the geometry and topology of
the family over the punctured disk: $f^*:\X^*=\X\setminus f^{-1}(0)
\rightarrow S^*=S\setminus\{0\}$. For example, in some sense, the
degree to which $\X_0$ is singular can be measured in terms of the
monodromy operator $T:H^*(\X_t,\ZZ)\rightarrow H^*(\X_t,\ZZ)$, where
$t\in S^*$ is a basepoint of a simple loop around the origin of $S$.

An appropriate form of the mirror symmetry conjecture
suggests that
associated to any sufficiently ``bad'' degeneration of Calabi-Yau
manifolds, i.e. a maximally unipotent degeneration or large complex
structure limit point, there should be a mirror manifold $\check X$,
defined as a symplectic manifold. Furthermore, if the family 
$f^*:\X^*\rightarrow S^*$ is polarized, i.e. given a choice
of a relatively ample line bundle $\shL$ on $\X^*$, one should expect to be able
to construct a degenerating family of complex manifolds $\check\X^*
\rightarrow S^*$ along with a polarization $\check\shL$. This correspondence
between degenerating polarized families is not precise, though
it can be made more precise if one allows multi-parameter families
of Calabi-Yau manifolds and subcones of the relatively ample
cone of $f^*:\X^*\rightarrow S^*$: this is essentially the form in 
which a general mirror symmetry conjecture was stated in \cite{Morlcsl}
and we do not wish to elaborate on that point of view here.

Instead, let us ask a number of questions which arise once one begins to
think about degenerations of Calabi-Yau manifolds:

(1) Mumford \cite{Mumf} constructed degenerations of $n$ complex dimensional
abelian varieties using combinatorial data consisting of periodic
polyhedral decompositions of $\RR^n$. This has been developed further
by Faltings and Chai \cite{FC}, Alexeev and Nakamura \cite{AN}, and 
Alexeev \cite{Alex}. Is there an analogue of this construction for
Calabi-Yau manifolds? What kind of combinatorial data would be
required to specify a degeneration?

(2) If one views a degeneration $\X\rightarrow S$ of Calabi-Yau manifolds
as a smoothing of a singular variety $\X_0$, can we view
mirror symmetry as an operation on {\it singular varieties}, 
exchanging $\X_0$ and $\check\X_0$, with the smoothings of these
singular varieties being mirror to each other?

(3) When one has a degeneration $f:\X\rightarrow S$, it is a standard
fact of life of birational geometry that there may be other birationally
equivalent families $f':\X'\rightarrow S$, with
$f^{-1}(S\setminus\{0\})\cong (f')^{-1}(S\setminus\{0\})$ but $\X\not\cong
\X'$. In other words, there may be many different ways of putting in
a singular fibre. It is traditional to take a semistable
degeneration, where $\X_0$ has normal crossings, but this
may not be the most natural thing to do in the context of
mirror symmetry. Is there a more natural class of compactifications
to consider?

(4) Can we clearly elucidate the connection between the singular fibre
$\X_0$ and the topology and geometry of the conjectural
Strominger-Yau-Zaslow fibration on $\X_t$ for $t$ small?

Let us look at a simple example to clarify these questions. Consider the
family of K3 surfaces $\X\subseteq \PP^3\times S$ given by
$tf_4+z_0z_1z_2z_3=0$, where $t$ is a coordinate on the disk $S$
and $z_0,\ldots,z_3$ are homogeneous coordinates on $\PP^3$. Here
$f_4$ is a general homogeneous polynomial of degree $4$. Then
$\X_0$ is just the union of coordinate planes in $\PP^3$,
and the total space of $\X$ is singular at the $24$ points of 
$\{t=f_4=0\}\cap Sing(\X_0)$. Thus, while $\X_0$ is normal crossings,
$\X\rightarrow S$ is not semi-stable as $\X$ is singular, and it
is traditional to obtain a semi-stable degeneration by blowing
up the irreducible components of $\X_0$ in some order. However,
this makes the fibre over $0$ much more complicated, and we
would prefer not to do this: before blowing up, the irreducible
components of $\X_0$ are toric varieties meeting along toric strata,
but after blowing up, the components become much more complex.

A second point is that we would like to allow $\X_0$ to have
worse singularities than simple normal crossings. A simple
example of why this would be natural is a generalisation of the
above example. Let $\Xi\subseteq\RR^n$ be a reflexive polytope,
defining a projective toric variety $\PP_{\Xi}$, along with a line
bundle $\O_{\PP_{\Xi}}(1)$. Let $s\in \Gamma(\PP_{\Xi},\O_{\PP_{\Xi}}(1))$
be a general section. There is also a special section $s_0$,
corresponding to the unique interior point of $\Xi$, and $s_0$
vanishes precisely on the toric boundary of the toric variety
$\PP_{\Xi}$, i.e. the complement of the big $(\CC^*)^n$ orbit in $\PP_{\Xi}$.
Then, as before, we can consider a family $\X\subseteq\PP_{\Xi}\times S$
given by the equation $ts+s_0=0$. Then $\X_0$ is the toric boundary
of $\PP_{\Xi}$, but is not necessarily normal crossings. It then
seems natural, at the very least, to allow $\X_0$ to locally look like
the toric boundary of a toric variety. Such an $\X_0$
is said to have {\it toroidal crossings singularities}. (\cite{SS2})

These two generalisations of normal crossings, i.e. allowing
the total space $\X$ to have some additional singularities and allowing
$\X_0$ to have toroidal crossings, will provide a natural category
of degenerations in which to work. In particular, we introduce
the notion of a {\it toric degeneration} of Calabi-Yau manifolds
in \S 1, which will formalise the essential features of the above
examples. A toric degeneration $f:\X\rightarrow S$ will, locally
away from some nice singular set $Z\subseteq\X$, look as if it is
given by a monomial in an affine toric variety. Furthermore, the
irreducible components of $\X_0$ will be toric varieties meeting
along toric strata. See Definition \ref{toricdegen} for a precise
definition.

Toric degenerations provide a natural generalisation of Mumford's
degenerations of abelian varieties. The key point for us is
that there is a generalisation of the combinatorics necessary to 
specify a degeneration of abelian varieties which can be applied
to describe toric degenerations of Calabi-Yau manifolds. To any toric
degeneration $f:\X\rightarrow S$ of $n$ complex dimensional Calabi-Yau
manifolds, one can build the dual intersection complex $B$ of $\X_0$, which
is a cell complex of real dimension $n$. If $\X_0$ were normal crossings,
this would be the traditional dual intersection graph, which is
a simplicial complex. The nice thing is that using the toric data
associated to the toric degeneration, we can put some additional data
on $B$, turning it into an {\it integral affine manifold with singularities}.
An integral affine manifold is a manifold with coordinate charts whose
transition functions are integral affine transformations. An integral
affine manifold with singularities will be, for us, a manifold
with an integral affine structure off of a codimension two subset.

Under this affine structure, the cells of $B$ will be polyhedra, and the
singularities of $B$ will be intimately related to the singular set
$Z\subseteq\X$. In the case of a degeneration of abelian varieties,
$B$ is just $\RR^n/\Lambda$, where $\Lambda\subseteq\ZZ^n$ is a lattice.
The polyhedral decomposition of $B$ is the same data required for
Mumford's construction. This combinatorial construction will
be explained in \S 2.

To complete the answer to (1), we need a way of going backwards: given
an integral affine manifold with singularities along with a
suitable decomposition into polyhedra, is it possible to construct
a degeneration $f:\X\rightarrow S$ from this? It turns
out to be easy to construct the singular fibre from this data.
To construct an actual degeneration, we need to use deformation
theory and try to smooth $\X_0$. However, this cannot be done
without some additional data on $\X_0$. Indeed, there may be many
distinct ways of smoothing $\X_0$. To solve this problem, one must
place a {\it logarithmic structure} of Illusie--Fontaine on 
$\X_0$. It is difficult to gain an intuition for logarithmic structures,
and we will try to avoid doing too much log geometry in this announcement. 
Suffice it
to say in this introduction that there is some additional structure
we can place on $\X_0$ turning it into a {\it log scheme}, which we
write as $\X_0^{\dagger}$. This preserves some of the information
associated to the inclusion $\X_0\subseteq\X$, but can be described
without knowing $\X$. It is then hoped, and definitely is the case
in two dimensions, that a log scheme $\X_0^{\dagger}$
has good deformation theory and can be deformed, and in some
cases yield a family $\X\rightarrow S$ as desired. This is
as yet the most technically difficult aspect of our program, and has yet to
be fully elaborated.

If the above ideas can be viewed as a generalisation of Mumford's
construction, they in fact give the key idea to answer question (2).
In fact, we cannot define a mirror symmetry operation
which works on singular fibres unless we incorporate a log
structure. One of the fundamental points of our program will be to define
mirrors of log Calabi-Yau schemes of the sort arising from toric
degenerations of Calabi-Yau manifolds.

The main point of this construction is as follows. Given a toric
degeneration $\X\rightarrow S$, or just as well $\X_0^{\dagger}$,
we obtain an affine manifold $B$ along with a polyhedral
decomposition $\P$. If $\X$, or $\X_0^{\dagger}$, is equipped with
an ample line bundle, we can also construct the intersection complex
of $\X_0$. While vertices of the dual intersection complex correspond
to irreducible components of $\X_0$, it is the maximal cells of the
intersection complex which correspond to these same components. However,
the presence of the polarization allows us to define an affine structure
with singularities on the intersection complex. Thus we have
a new affine manifold
$\check B$, and a new polyhedral decomposition $\check\P$. Most significantly,
there is a clear relationship between $B,\P$ and $\check B,\check\P$.
The polarization on $\X_0^{\dagger}$ can be viewed as giving a convex
multi-valued piecewise linear function $\varphi$ on $B$, and one can define the
{\it discrete Legendre transform} of the triple $(B,\P,\varphi)$, being
a triple $(\check B,\check\P,\check\varphi)$, where $\check B,\check\P$
are as above and $\check\varphi$ is a convex multivalued piecewise linear
function on $\check B$. These affine structures are naturally dual.
This discrete Legendre transform is analagous to the well-known
discrete Legendre transform of a convex piecewise linear function on
$\RR^n$. 

From the data $(\check B,\check\P,\check\varphi)$ one can then construct
a scheme $\check\X_0$ along with a 
log structure. If $\X_0^{\dagger}$ and $\check\X_0^{\dagger}$ are
smoothable, then we expect that the smoothings will lie in mirror
families of Calabi-Yau manifolds.

This construction can be seen to reproduce the Batyrev-Borisov
mirror symmetry construction for complete intersections in toric
varieties. This is evidence that it is the correct construction.
More intriguingly, in a strong sense this is an
algebro-geometric analogue of the Strominger-Yau-Zaslow conjecture.

Indeed, we show that given a toric degeneration $\X\rightarrow S$, if
$B$ is the corresponding integral affine manifold with singularities,
then the integral affine structure on $B$ determines a torus bundle
over an open subset of $B$, and topologically $\X_t$ is
a compactification of this torus bundle. Furthermore, it is well-known
(\cite{Hit}, \cite{Leung}) that such affine structures can be dualised
via a (continuous) Legendre transform: we are just replacing the continuous
Legendre transform with a discrete Legendre transform. As a result, our
construction gives both a proof of a topological form of
SYZ in a quite general context, including the Batyrev-Borisov complete
intersection case, and also points the way towards an algebro-geometrization
of the SYZ conjecture.

This paper is meant to serve as an announcement of these ideas, which still
are a work in progress. We will
give an outline of the approach, and suggest what may be proved using
it. A full exploration of this program is currently ongoing. Details
of many of the ideas discussed here will appear in \cite{GS}.

Some aspects of the ideas here were present in earlier literature. The
idea that mirror symmetry can be represented by an exchange of
irreducible components and information about deepest points was present
in an intuitive manner in Leung and Vafa's paper \cite{LV}. The idea of using 
the dual intersection complex (in the normal crossing case) to represent
the base of the SYZ fibration first occurred in Kontsevich and Soibelman's
paper \cite{KS}, and is expanded on in Kontsevich's ideas of
using Berkovich spaces \cite{K}. 
Finally, the idea that logarithmic structures can play
an important role in mirror symmetry first appears in \cite{SS2}, which
served as a catalyst for this work. 

{\it Acknowledgements:} We would like to thank Maxim Kontsevich, Simone
Pavanelli, Richard Thomas, Stefan Schr\"oer and Ilia Zharkov
for useful discussions.

\section{Toric degenerations}

\begin{definition}
\label{toricdegen}
Let $f:\X\rightarrow S$ be a proper flat family of relative dimension
$n$, where $S$ is a disk and $\X$ is a complex analytic space
(not necessarily non-singular). We say $f$ is a {\it toric degeneration}
of Calabi-Yau varieties if 
\begin{enumerate}
\item the canonical bundle $\omega_{\X}$ of $\X$ is trivial.
\item
$\X_t$ is an irreducible normal Calabi-Yau variety with only
canonical singularities for $t\not=0$. (The reader may
like to assume $\X_t$ is smooth for $t\not=0$).
\item $\X_0$ is a reduced Cohen-Macaulay variety, with
normalisation $\nu:\tilde\X_0\rightarrow\X_0$. In addition, $\tilde\X_0$
is a disjoint union of toric varieties, and if $U\subseteq\tilde\X_0$
is the union of big $(\CC^*)^n$ orbits in $\tilde\X_0$, then $\nu:U
\rightarrow \nu(U)$ is an isomorphism and $\nu^{-1}(\nu(U))=U$. Thus $\nu$ only
identifies toric strata of $\tilde\X_0$. Furthermore, if
$S\subseteq \tilde\X_0$ is a toric stratum, then $\nu:S\rightarrow \nu(S)$
is the normalization of $\nu(S)$.
\item There exists a closed subset $Z\subseteq\X$ of relative 
codimension $\ge 2$ such that $Z$ satisfies the following properties:
$Z$ does not contain the image under $\nu$ 
of any toric stratum of $\tilde\X_0$, 
and for any point $x\in \X\setminus Z$, there is a neighbourhood
$U_x$ (in the analytic topology) of $x$, an
$n+1$-dimensional affine toric variety $Y_x$, a regular function
$f_x$ on $Y_x$ given by a monomial, and a commutative diagram
$$\begin{matrix}
U_x&\mapright{\psi_x}&Y_x\cr
\mapdown{f|_{U_{x}}}&&\mapdown{f_x}\cr
S&\mapright{\varphi_x}&\CC\cr
\end{matrix}$$
where $\psi_x$ and $\varphi_x$ are open embeddings. Furthermore,
$f_x$ vanishes precisely once on each toric divisor of $Y_x$.
\end{enumerate}
\end{definition}

There are two key toric aspects needed here: each irreducible component
(or normalization thereof) of $\X_0$ is toric meeting other
components only along toric strata, and $f$ on a
neighbourhood of each point away from $Z$ looks like a morphism from
a toric variety to the affine line given by a monomial of a special sort. 
It is important
to have both conditions, as mirror symmetry will actually exchange
these two bits of toric data.

Key examples to keep in mind were already given in the introduction:
a degeneration of quartics in $\PP^3$ to $z_0z_1z_2z_3=0$, or a family
of hypersurfaces in a toric variety, degenerating to the variety
$s_0=0$, given by the equation $st+s_0=0$. As long as $s$
does not vanish on any toric stratum of $\PP_{\Xi}$, then the degeneration
given is toric, with the singular set $Z$ given by
$$Z=\X\cap [(Sing(\X_0)\cap \{s=0\})\times S]\subseteq \PP_{\Xi}\times S.$$
We will keep these examples in mind as we continue through the paper.

There are more general forms of these degenerations also: instead of
weighting a general section $s$ with a single factor $t$,
one can consider equations of the form
$$s_0+\sum_{\hbox{$s\in\Xi$ integral}} t^{h(s)+1}s=0$$
where $h$ is a suitable height function on the set of integral points
of $\Xi$. Here we have identified the set of integral points of $\Xi$
with a monomial basis for $\Gamma(\PP_{\Xi},\O_{\PP_{\Xi}}(1))$ as usual.
Now the variety $\X$ defined by this equation in $\PP_{\Xi}\times S$
is in general too singular to give rise to a toric degeneration;
however, there are standard techniques for obtaining a partial
desingularization of $\X$ to yield toric degenerations. This sort
of technique has been applied in \cite{Zhar} and \cite{HZ}, and a variant of
this in \cite{Hu}, where the goal was to obtain semi-stable degenerations
instead of toric degenerations.

There are other examples of toric degenerations: any maximally
unipotent degeneration of abelian varieties is toric, and examples
of toric degenerations of Kodaira surfaces are given in \cite{SS}.

\section{From toric degenerations to affine manifolds: the dual intersection
complex}

We fix
$M=\ZZ^n$, $N=\Hom_{\ZZ}(M,\ZZ)$, $M_{\RR}=M\otimes_{\ZZ}
\RR$, $N_{\RR}=N\otimes_{\ZZ} {\RR}$. We set
$$\Aff(M_{\RR})=M_{\RR}\rtimes GL_n(\RR)$$
to be the group of affine transformations of $M_{\RR}$, with subgroups
\begin{eqnarray*}
\Aff_{\RR}(M)&=&M_{\RR}\rtimes GL_n(\ZZ)\\
\Aff(M)&=&M\rtimes GL_n(\ZZ).
\end{eqnarray*}

\begin{definition}
Let $B$ be an $n$-dimensional manifold.
An {\it affine structure} on $B$ is given by an open cover $\{U_i\}$
along with coordinate charts $\psi_i:U_i\rightarrow M_{\RR}$,
whose transition functions $\psi_i\circ\psi_j^{-1}$ lie in $\Aff(M_{\RR})$.
The affine structure is {\it integral} if the transition functions
lie in $\Aff(M)$. 
If $B$ and $B'$ are (integral) affine
manifolds of dimension $n$ and $n'$ respectively,
then a continuous map $f:B\rightarrow B'$ is (integral) affine if
locally $f$ is given by affine linear transformations from $\RR^n$ to
$\RR^{n'}$ ($\ZZ^n$ to 
$\ZZ^{n'}$). 
\end{definition}

\begin{definition}
 An affine manifold with singularities is a $C^0$
(topological) manifold $B$ along with a set $\Delta\subseteq B$ which
is a finite union of locally closed submanifolds of codimension at
least 2, and an affine structure on $B_0=B\setminus \Delta$.
An affine manifold with singularities is {\it integral}
if the affine structure on $B_0$ is integral. 
We always denote by $i:B_0\hookrightarrow B$ the inclusion map.
A continuous
map $f:B\rightarrow B'$ of (integral) affine manifolds with singularities is
(integral) affine  if $f^{-1}(B_0')\cap B_0$ is dense in $B$ and
$$f|_{f^{-1}(B_0')\cap B_0}:f^{-1}(B_0')\cap B_0\rightarrow B_0'$$
is (integral) affine.
\end{definition}

\begin{definition} 
\label{PD1}
A {\it polyhedral decomposition} of a closed set $R\subseteq
M_{\RR}$ is a locally finite covering $\P$ of $R$ by closed convex
polytopes (called {\it cells}) with the property that
\item{(1)} if $\sigma\in\P$ and $\tau\subseteq\sigma$ is a face of
$\sigma$ then
$\tau\in\P$;
\item{(2)} if $\sigma,\sigma'\in\P$, then $\sigma\cap\sigma'$ is a common
face of $\sigma$ and $\sigma'$.

We say the decomposition is {\it integral} if all vertices (0-dimensional
elements of $\P$) are contained in $M$.
\end{definition}

For a polyhedral decomposition $\P$ and $\sigma\in\P$ we define
$$Int(\sigma)=\sigma\setminus \bigcup_{\tau\in\P,\tau\subsetneq\sigma}\tau.$$

We wish to define a polyhedral decomposition of an affine manifold
with singularities generalizing the above notion for a set in $M_{\RR}$.
This will be a decomposition of $B$ into lattice polytopes with respect
to the integral affine structure on $B$. This definition must be
phrased rather carefully, as we need to control the interaction
between these polytopes and the discriminant locus $\Delta$ of $B$.
In particular, $\Delta$ should contain no zero-dimensional cells
and not pass through the interior of any $n$-dimensional cell, but
there are subtler restrictions necessary for our purposes.
We also need to allow in general for cells to have self intersection.
For example, by identifying opposite sides as depicted,
the following picture shows a polyhedral decomposition
of $B=\RR^2/\ZZ^2$ with one two-dimensional cell, two one-dimensional
cells, and one zero-dimensional cell.
$$\epsfbox{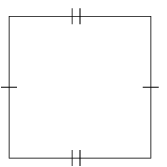}$$

Our definition is

\begin{definition}
\label{PD2}
Let $B$ be an integral affine manifold with singularities. A {\it polyhedral
decomposition} of $B$ is a locally finite covering $\P$ of $B$
by closed subsets of $B$ (called {\it cells})
which satisfies the following properties. If $\{v\}\in\P$ for some point
$v\in B$, then $v\not\in\Delta$ and
there exists an integral polyhedral decomposition $\P_v$
of a closed set $R_v\subseteq M_{\RR}$ which is the closure of an open
neighbourhood of the origin, and a
continuous map ${\rm exp}_v:R_v\rightarrow B$, ${\rm exp}_v(0)=v$, satisfying
\item{(1)} ${\rm exp}_v$ 
is locally an immersion (of manifolds with boundary) onto its image and
is an integral affine map in some neighbourhood of the origin. 
\item{(2)} For every $n$-dimensional $\tilde\sigma\in\P_v$, 
${\rm exp}_v(Int(\tilde\sigma))\cap\Delta=\phi$ and the restriction
of ${\rm exp}_v$ to $Int(\tilde\sigma)$ is integral affine. 
\item{(3)} $\hbox{$\sigma\in\P$ and $v\in \sigma$}
\Leftrightarrow
\hbox{$\sigma={\rm exp}_v(\tilde\sigma)$ for some $\tilde\sigma\in\P_v$ with
$0\in\tilde\sigma$.}$
\item{(4)} Every $\sigma\in\P$ contains a point $v\in\sigma$ with 
$\{v\}\in\P$.

In addition we say the polyhedral decomposition is {\it toric} if it satisfies
the additional condition
\item{(5)}
For each $\sigma\in\P$, there is a neighbourhood $U_{\sigma}
\subseteq B$
of $Int(\sigma)$ and an integral
affine submersion $s_{\sigma}:U_{\sigma}\rightarrow M'_{\RR}$ where
$M'$ is a lattice of rank equal to
$\dim B-\dim \sigma$ and $s_{\sigma}(\sigma\cap
U_{\sigma})=\{0\}$.
\end{definition}

\begin{example} 
\label{torusexample}
If $B=M_{\RR  }$, $\Delta=\phi$,
a polyhedral decomposition of $B$ is just an integral polyhedral
decomposition of $M_{\RR}$ in the sense of Definition
\ref{PD1}. If $B=M_{\RR}/\Lambda$ for some lattice $\Lambda\subseteq M$, 
then a polyhedral decomposition of $B$ is induced by
a polyhedral decomposition of $M_{\RR}$ invariant under $\Lambda$.

Explicitly, in the example given above for $B=\RR^2/\ZZ^2$, we have a
unique vertex $v$, and we can take $R_v$ to be a union of four copies
of the square in $\RR^2$:
$$\epsfbox{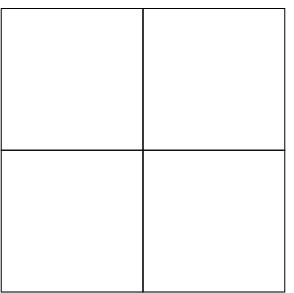}$$
We take ${\rm exp}_v$ to be the restriction to $R_v$ of the quotient
map $\RR^2\rightarrow\RR^2/\ZZ^2$. Of course, we could have taken
$R_v=M_{\RR}$ also in this example. Note that in this example,
${\rm exp}_v$ is not even an isomorphism in the interior of $R_v$
but just an immersion.

In the case where $\Delta$ is empty, the toric condition
is vacuous.
When $B$ has singularities, the definition of toric polyhedral decomposition
imposes some slightly subtle restrictions on how the cells of $\P$
interact with $\Delta$. 
\end{example}

\begin{remark}
\label{fan1}
Given a polyhedral decomposition $\P$ on $B$, if $v$ is a vertex of $\P$,
we can look at the polyhedral decomposition of $R_v$ in a small neighbourhood
of the origin in $M_{\RR}$. This clearly coincides with
a small neighbourhood of the origin of a complete rational polyhedral
fan $\Sigma_v$ in $M_{\RR}$:
$$\epsfbox{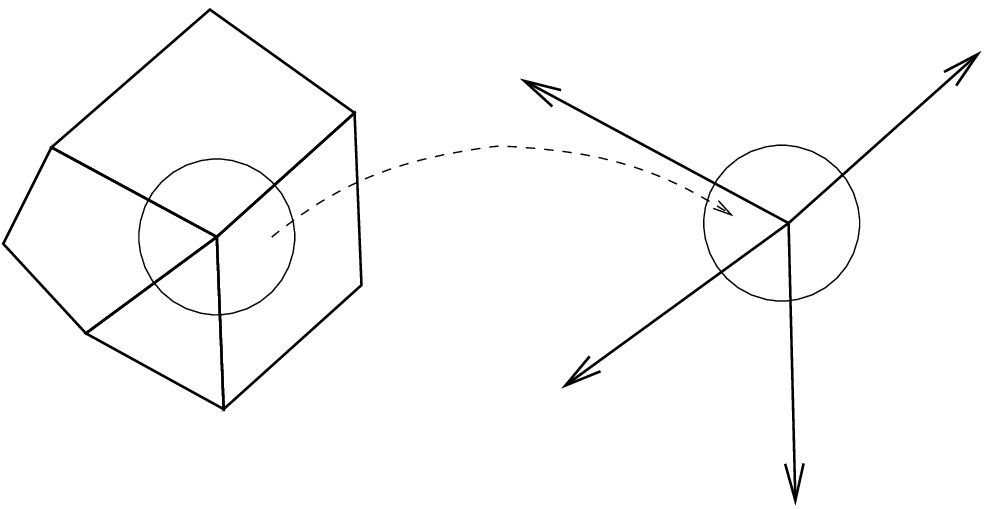}$$
In fact, we shall now see that the data of the affine structures on
maximal cells of $\P$ and this fan structure at each vertex
$v$ essentially determine the affine structure on $B$.
\end{remark}

\begin{definition} 
Recall that if $\sigma\subset M_{\RR}$ is a polytope, then the barycentre 
$Bar(\sigma)$ of $\sigma$ is the average of the vertices of $\sigma$, and thus
is invariant under affine transformations. The first barycentric subdivision
of $\sigma$ is then the triangulation of $\sigma$ consisting of
all simplices spanned by barycentres of ascending chains of cells of $\sigma$.
Thus given a polyhedral decomposition $\P$ of an affine manifold with
singularities $B$, we can define the first barycentric subdivision
$Bar(\P)$ of $\P$ to be the triangulation
consisting of all images of simplices in the first barycentric subdivisions
of all $\tilde\sigma\in \P_v$ for all vertices $v$. Because barycentric
subdivisions are affine invariants, this gives a well-defined triangulation
of $B$.
\end{definition}

We will now describe a standard procedure for constructing affine
manifolds with singularities along with polyhedral decompositions.
Let $\P'$ be a collection of $n$-dimensional
integral polytopes in $M_{\RR}$. Suppose we are given integral
affine identifications of various faces of the polytopes
in $\P'$ in such a way that once we glue the polytopes using these
identifications, we obtain a manifold $B$, along with a decomposition
$\P$ consisting of images of faces of polytopes in $\P'$. In particular,
we have the identification map
$$\pi:\coprod_{\sigma'\in\P'} \sigma'\rightarrow B.$$

Now $B$ is not yet an affine manifold with singularities. It only
has an affine structure defined in the interiors of maximal cells.
When two polytopes of $\P'$ are identified along subfaces, we have
an affine structure on that subface, but no affine structure in
the directions ``transversal'' to that subface. We cannot,
however, expect an affine structure on all of $B$, and we need to choose
a discriminant locus. We do this as follows.

Let $Bar(\P)$ be the first barycentric subdivision of $\P$.
We define $\Delta'\subseteq B$ to be the union of all simplices 
in $Bar(\P)$ not containing a vertex of $\P$ or the barycentre of
an $n$-dimensional cell of $\P$. This can be seen as the codimension
two skeleton of the dual cell complex to $\P$.

For a vertex $v$ of $\P$, let $W_v$ be the union of the interiors
of all simplices in $Bar(\P)$ containing $v$. Then $W_v$ is an open
neighbourhood of $v$, and
$$\{W_v|\hbox{$v$ a vertex of $\P$}\}\cup \{Int(\sigma)|
\hbox{$\sigma$ a maximal cell of $\P$}\}$$
form an open covering of $B\setminus\Delta'$. To define an affine structure
on $B\setminus\Delta'$, we need to choose affine charts on $W_v$.

For a vertex $v$ of $\P$, let
$$\P'_v=\{(v',\sigma')|\hbox{$v'\in\sigma'\in\P'$ a vertex, $\pi(v')=v$}\}.$$
Let $R_v$ be the quotient of $\coprod_{(v',\sigma')\in\P'_v}\sigma'$
by the equivalence relation which identifies proper faces 
$\omega_i'\subsetneq\sigma_i'$
for $i=1,2$, $(v_i',\sigma_i')\in\P'_v$,
if $\pi(\omega_1')=\pi(\omega_2')$ and $v_i'\in\omega_i'$. For example, if
$\P'$ consists of the unit square in $\RR^2$, and $B$ is obtained by identifying
opposite sides, we have a unique vertex $v$ in $\P$ and the picture
$$\epsfbox{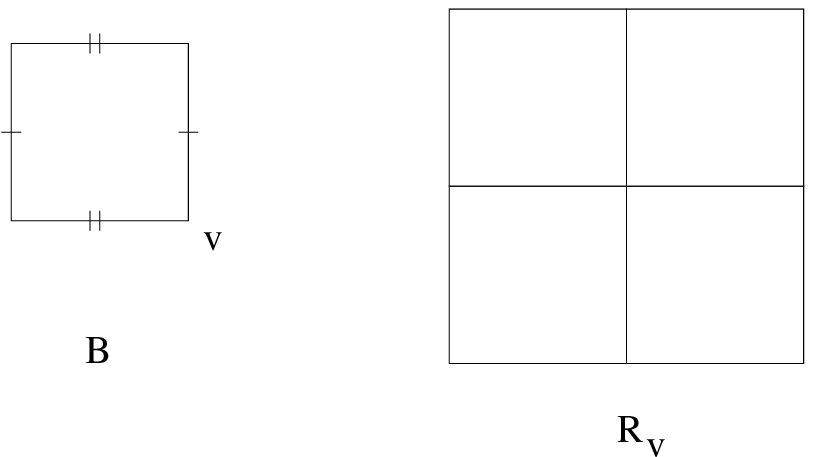}$$
as in Example \ref{torusexample}
There is a continuous map 
$$\pi_v:R_v\rightarrow B$$ defined by taking $b\in\sigma'\subseteq R_v$
to $\pi(b)$, and it is easy to see that if $U_v$ is
the connected component of $\pi_v^{-1}(W_v)$ containing $0$,
then $U_v\rightarrow W_v$ is
a homeomorphism.

$R_v$ has an abstract polyhedral decomposition $\P_v$, with a unique
vertex $v'\in U_v$ mapping to $v$. We will need to find an
embedding $i_v:R_v\rightarrow M_{\RR}$. If this is done in the appropriate
way, then a coordinate chart $\psi_v:W_v\rightarrow M_{\RR}$ can
be defined as $i_v\circ\pi_v^{-1}|_{W_v}$, and ${\rm exp}_v:i_v(R_v)\rightarrow
B$ can be defined as $\pi_v\circ i_v^{-1}$, giving both an affine structure
on $B\setminus\Delta'$ and a proof that $\P$ is a polyhedral decomposition
of $B$.

To do this, we need to choose a {\it fan structure} at each vertex $v$
of $\P$. This means for each $v$ we choose a complete rational polyhedral
fan $\Sigma_v$ in $M_{\RR}$ and a one-to-one inclusion preserving
correspondence between elements of $\P_v$ containing $v'$ and
elements of $\Sigma_v$ which we write as $\sigma\mapsto \sigma_{v'}$. 
Furthermore, this correspondence should have the property that
there exists an integral affine isomorphism $i_{\sigma}$ between the tangent
wedge of $\sigma$ at $v'$ and $\sigma_{v'}$ which preserves the correspondence.
Such an isomorphism, if it exists, is unique (integrality is essential
here as otherwise we can rescale). By this uniqueness, the
maps $i_{\sigma}$ glue together to give a map
$$i_v:R_v\rightarrow M_{\RR}$$ which is a homeomorphism onto its image.
Then it is easy to see that using $\psi_v$ and ${\rm exp}_v$ as defined above
one obtains an integral affine structure on $B_0$ and one sees that $\P$
is a polyhedral decomposition.

Suppose we have made such choices of fan structure, 
and so obtained an affine structure
on $B\setminus\Delta'$ and a polyhedral decomposition $\P$. It often happens
that our choice of $\Delta'$ is too crude, and we can still extend the
affine structure to a larger open set of $B$. 
Let $\Delta$ be the smallest subset of $\Delta'$ such that the affine
structure on $B\setminus\Delta'$ extends to $B\setminus\Delta$
(Such an extension is unique if it exists).
The set $\Delta$ can be characterized precisely, but we will not do
this here. We call $\Delta$ the {\it minimal discriminant locus}.
\bigskip

The main point of this section is that given a toric degeneration
of Calabi-Yau manifolds $f:\X\rightarrow S$, there is a natural integral
affine manifold with singularities $B$ we can associate to it, the dual 
intersection complex.

We will construct $B$ as a union of lattice polytopes as in \S 2,
specifying a fan structure at each vertex. Specifically,
let the normalisation of $\X_0$, $\tilde \X_0$, be written as a disjoint
union $\coprod X_i$ of toric varieties $X_i$, $\nu:\tilde\X_0\rightarrow\X_0$
the normalisation. The {\it strata} of $\X_0$ are the elements of
the set 
$$Strata(\X_0)=\{\nu(S)|\hbox{$S$ is a toric stratum of $X_i$ for some $i$}\}.$$
Here by toric stratum we mean the closure of a $(\CC^*)^n$ orbit.

Let $\{x\}\in Strata(\X_0)$ be a zero-dimensional stratum. Let
$M'=\ZZ^{n+1}$, $M'_{\RR}=M'\otimes_{\ZZ} \RR$, $N'=\Hom_{\ZZ}(M',\ZZ)$,
$N'_{\RR}=N'\otimes_{\ZZ}\RR$ as usual. Then applying 
Definition \ref{toricdegen} (4) to a neighbourhood of $x$, 
there is a toric variety $Y_x$ defined by a rational polyhedral
cone $\tilde\sigma_x\subseteq M'_{\RR}$ such that in
a neighbourhood of $x$, $f:\X\rightarrow S$ is locally isomorphic to
$f_x:Y_x\rightarrow\CC$, where $f_x$ is given by a monomial, i.e.
an element $\rho_x\in N'$. Now put
$$\sigma_x=\{m\in\tilde\sigma_x| \langle \rho_x,m\rangle=1\}.$$
Recall that there is a one-to-one correspondence between codimension
one toric strata of $Y_x$ and the dimension one faces of $\tilde\sigma_x$:
these strata are precisely the toric divisors of $Y_x$. Now the condition
that $f_x$ vanishes to order 1 on each such divisor can be expressed
as follows. For every one-dimensional face $\tilde\tau$ of
$\tilde\sigma_x$, let $\tau$ be a primitive integral generator of
$\tilde\tau$. Then the order of vanishing of $f_x$ on the toric divisor
corresponding to $\tilde\tau$ is $\langle\rho_x,\tau\rangle$. Since this
must be 
$1$, we see in fact that $\sigma_x$ is the convex hull of the primitive
integral generators of the one-dimensional faces of $\tilde\sigma_x$.
If we put $M=\{m\in M'|\langle \rho_x,m\rangle=1\}$,
$M_{\RR}=\{m\in M'_{\RR}|\langle\rho_x,m\rangle=1\}$, then $\sigma_x$
is a lattice polytope in the affine space $M_{\RR}$.

What is less obvious, but which follows from the triviality of the
canonical bundle of $\X_0$, is that $\sigma_x$ is in fact an $n$-dimensional
lattice polytope.

\begin{example} If at a point $x\in \X_0$ which is a zero-dimensional
stratum, the map $f:\X\rightarrow S$ is locally isomorphic to 
$\CC^{n+1}\rightarrow \CC$ given by $(z_0,\ldots,z_n)\mapsto \prod_{i=0}^n z_i$,
we say $f$ is normal crossings at $x$. In this case, the relevant
toric data is as follows: $\tilde\sigma_x$ is generated by the points
$(1,0,\ldots,0),\ldots,(0,\ldots,0,1)$ in $\RR^{n+1}$, and the map
is given by the monomial determined by $(1,1,\ldots,1)\in\Hom_{\ZZ}(\ZZ^{n+1},
\ZZ)$. Then $\sigma$ is the standard simplex in the affine space
$$\{(x_0,\ldots,x_n)|\sum_{i=0}^n x_i=1\}$$
with vertices the standard basis of $\RR^{n+1}$.
\end{example}

We can now describe how to construct $B$ by gluing together the polytopes
$$\{\sigma_x| \{x\}\in Strata(\X_0)\}.$$
We will do this in the case that every irreducible component
of $\X_0$ is in fact itself normal so that $\nu:X_i\rightarrow \nu(X_i)$ is an
isomorphism. The reader may be able to imagine the more
general construction.

Note in this case there is a one-to-one inclusion reversing
correspondence between faces of
$\sigma_x$ and elements of $Strata(\X_0)$ containing $x$. We can then
identify faces of $\sigma_x$ and $\sigma_{x'}$ if they correspond
to the same strata of $\X_0$. Some argument is necessary to show that this
identification can be done via an integral affine transformation, but
again this is not difficult.

Making these identifications, one obtains $B$. One can then prove

\begin{lemma} If $\X_0$ is $n$ complex dimensional, then $B$ is an $n$
real dimensional manifold.
\end{lemma}

The key point again here is the triviality of the canonical bundle.

As explained above, to give $B$ the structure of an affine manifold with 
singularities, we just need to specify a fan structure at each
vertex of $B$. Now the vertices $v_1,\ldots,v_m$ of $B$ are in 
one-to-one correspondence with the irreducible components $X_1,\ldots,X_m$
of $\X_0$ by construction. Each $X_i$ is toric, hence defined by
a complete fan $\Sigma_i$ living in $M_{\RR}$. For a vertex $v_i$,
the maximal cones of $\Sigma_i$ are in one-to-one correspondence with
the zero-dimensional strata of $X_i$. In fact, if $\sigma$
is a maximal cell of $B$ corresponding to a zero-dimensional
stratum of $X_i$, then $v_i\in\sigma$ and there is
a natural integral affine isomorphism between the corresponding
cone of $\Sigma_i$ and the tangent wedge of $\sigma$ at $v_i$.
It is not hard to see
this collection of isomorphisms gives a fan structure at each vertex $v_i$, thus
getting an integral affine manifold with singularities, along with 
a polyhedral decomposition $\P$.

\begin{examples}
\label{Bexamp}
(1)
In case $f$ is normal crossings away from the singular set $Z$, 
$B$ is the
traditional dual intersection complex: $B$ is a simplicial complex
with vertices $v_1,\ldots,v_m$
corresponding to the irreducible components
$X_1,\ldots,X_m$ of $\X_0$, and with a $p$-simplex with vertices $v_{i_0},
\ldots,v_{i_p}$ if $X_{i_0}\cap\cdots\cap X_{i_p}\not=\phi$. However,
the affine structure carries more information than the traditional
dual intersection complex because of the fan structure.

(2) Let $f:\X\rightarrow S$ be a degeneration of elliptic curves,
with $\X_0$ being a fibre of Kodaira type $I_m$, i.e. a cycle of $m$
rational curves. Furthermore, assume the total space $\X$ is 
non-singular. To ensure the irreducible components of
$\X_0$ are normal, we take $m\ge 2$. Then $f$ is normal crossings,
and $B$ is a cycle of $m$ line segments of length 1. There is a unique
fan structure here, with a neighbourhood of each vertex identified with
a neighbourhood of $0$ in the unique fan defining $\PP^1$:
$$\epsfbox{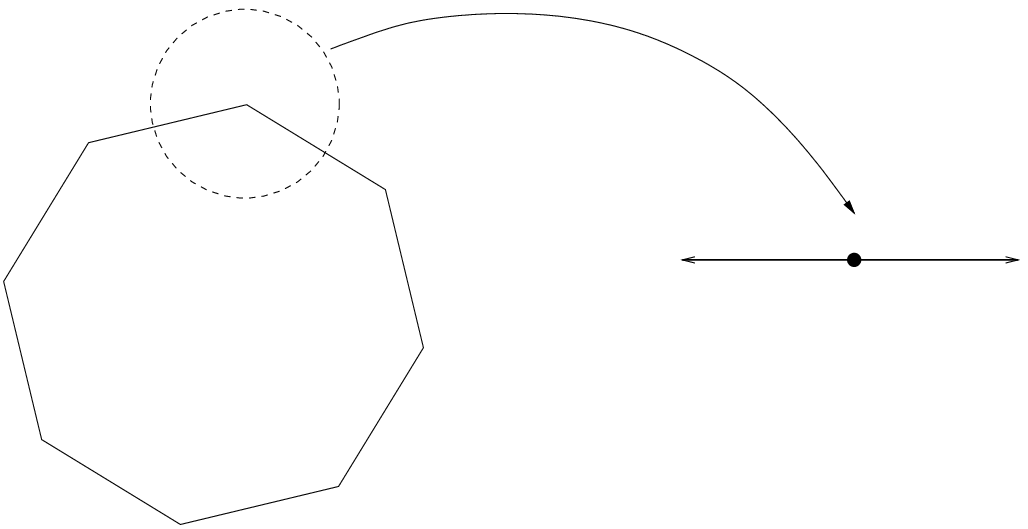}$$
As an affine manifold, $B$ is just $\RR/m\ZZ$, with the affine structure
induced by the standard one on $\RR$.

(3) To get a simple example which is not normal crossings,
one can start with the above example and contract some
chains of rational curves, so that $\X_0$ is still
a cycle of rational curves, but the total space $\X$ now has
singularities given locally by the equation $xy=z^n$ in $\CC^3$ for
various $n$ (where $n-1$ is the length of the chain contracted to create
the singularity). Locally, the map is given by $(x,y,z)\mapsto z$.
At such a point, the cone giving such a local description is
$\tilde\sigma$ generated by $(1,0)$ and $(1,n)$ in $\RR^2=M'_{\RR}$,
with the map given by $(1,0)\in N'$. Thus the corresponding 
$\sigma$ is a line segment of length $n$.
\begin{center}
\begin{picture}(0,0)%
\includegraphics{Ansing.pstex}%
\end{picture}%
\setlength{\unitlength}{1973sp}%
\begingroup\makeatletter\ifx\SetFigFont\undefined%
\gdef\SetFigFont#1#2#3#4#5{%
  \reset@font\fontsize{#1}{#2pt}%
  \fontfamily{#3}\fontseries{#4}\fontshape{#5}%
  \selectfont}%
\fi\endgroup%
\begin{picture}(3756,3928)(2923,-5084)
\put(3376,-1351){\makebox(0,0)[lb]{\smash{\SetFigFont{6}{7.2}{\rmdefault}{\mddefault}{\updefault}$(1,n)$}}}
\put(6316,-5026){\makebox(0,0)[lb]{\smash{\SetFigFont{6}{7.2}{\rmdefault}{\mddefault}{\updefault}$\sigma$}}}
\put(3361,-4921){\makebox(0,0)[lb]{\smash{\SetFigFont{6}{7.2}{\rmdefault}{\mddefault}{\updefault}$(1,0)$}}}
\end{picture}

\end{center}
Thus contracting a chain of $\PP^1$'s in (2) above has the effect of
keeping $B$ fixed but changing the polyhedral decomposition by erasing all 
vertices corresponding to these rational curves which have been contracted.

(4) A degenerating family of K3 surfaces: take $tf_4+x_0x_1x_2x_3=0$
in $\CC\times\PP^3$ as usual. Then $f$ is normal crossings at
each triple point of $\X_0$, so $B$ is obtained by gluing together
four standard simplices to form a tetrahedron. The chart for the
affine structure in a neighbourhood of a vertex $v$ identifies
that neighbourhood with a neighbourhood of zero of the fan
$\Sigma$ defining $\PP^2$; given the combinatorial correspondence
between the cells of $B$ containing $v$ and the cones of the fan
$\Sigma$, there is a unique such chart which is integral affine
on the interior of each $2$-cell containing $v$.
$$\epsfbox{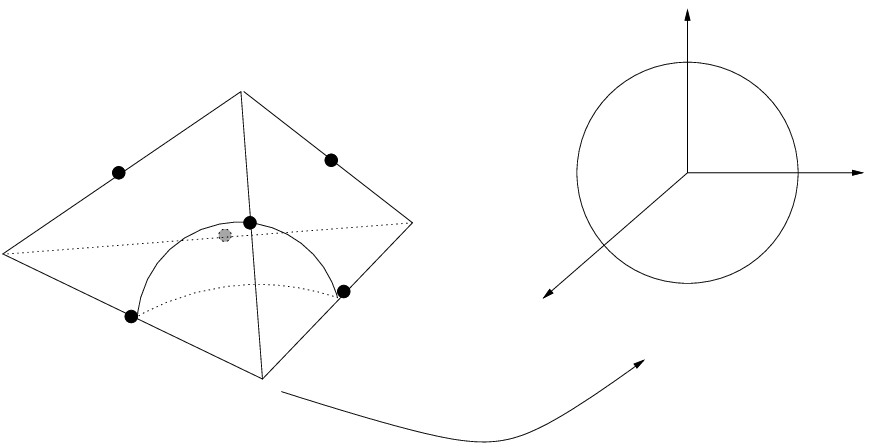}$$

In this example we have one singular point along each edge of the
tetrahedron indicated by the dots. 
These singularities cannot be removed. To see this, consider
two vertices $v_1$ and $v_2$. We can take $R_{v_1}$ and $R_{v_2}$
in the definition of polyhedral decomposition to look like
\begin{center}
\begin{picture}(0,0)%
\includegraphics{Rv12.pstex}%
\end{picture}%
\setlength{\unitlength}{2763sp}%
\begingroup\makeatletter\ifx\SetFigFont\undefined%
\gdef\SetFigFont#1#2#3#4#5{%
  \reset@font\fontsize{#1}{#2pt}%
  \fontfamily{#3}\fontseries{#4}\fontshape{#5}%
  \selectfont}%
\fi\endgroup%
\begin{picture}(8565,4144)(1111,-5429)
\put(2386,-2656){\makebox(0,0)[lb]{\smash{\SetFigFont{8}{9.6}{\rmdefault}{\mddefault}{\updefault}{\color[rgb]{0,0,0}$v_1$}%
}}}
\put(2521,-3016){\makebox(0,0)[lb]{\smash{\SetFigFont{8}{9.6}{\rmdefault}{\mddefault}{\updefault}{\color[rgb]{0,0,0}$(0,0)$}%
}}}
\put(3721,-2776){\makebox(0,0)[lb]{\smash{\SetFigFont{8}{9.6}{\rmdefault}{\mddefault}{\updefault}{\color[rgb]{0,0,0}$(1,0)$}%
}}}
\put(1111,-1441){\makebox(0,0)[lb]{\smash{\SetFigFont{8}{9.6}{\rmdefault}{\mddefault}{\updefault}{\color[rgb]{0,0,0}$(-1,1)$}%
}}}
\put(2686,-3331){\makebox(0,0)[lb]{\smash{\SetFigFont{8}{9.6}{\rmdefault}{\mddefault}{\updefault}{\color[rgb]{0,0,0}$\sigma_1$}%
}}}
\put(2056,-3181){\makebox(0,0)[lb]{\smash{\SetFigFont{8}{9.6}{\rmdefault}{\mddefault}{\updefault}{\color[rgb]{0,0,0}$\sigma_2$}%
}}}
\put(8221,-1456){\makebox(0,0)[lb]{\smash{\SetFigFont{8}{9.6}{\rmdefault}{\mddefault}{\updefault}{\color[rgb]{0,0,0}$v_1$}%
}}}
\put(8446,-2791){\makebox(0,0)[lb]{\smash{\SetFigFont{8}{9.6}{\rmdefault}{\mddefault}{\updefault}{\color[rgb]{0,0,0}$v_2$}%
}}}
\put(9676,-1501){\makebox(0,0)[lb]{\smash{\SetFigFont{8}{9.6}{\rmdefault}{\mddefault}{\updefault}{\color[rgb]{0,0,0}$(1,1)$}%
}}}
\put(8446,-1756){\makebox(0,0)[lb]{\smash{\SetFigFont{8}{9.6}{\rmdefault}{\mddefault}{\updefault}{\color[rgb]{0,0,0}$(0,1)$}%
}}}
\put(7096,-5371){\makebox(0,0)[lb]{\smash{\SetFigFont{8}{9.6}{\rmdefault}{\mddefault}{\updefault}{\color[rgb]{0,0,0}$(-1,-2)$}%
}}}
\put(2401,-4186){\makebox(0,0)[lb]{\smash{\SetFigFont{8}{9.6}{\rmdefault}{\mddefault}{\updefault}{\color[rgb]{0,0,0}$v_2$}%
}}}
\put(2386,-4501){\makebox(0,0)[lb]{\smash{\SetFigFont{8}{9.6}{\rmdefault}{\mddefault}{\updefault}{\color[rgb]{0,0,0}$(0,-1)$}%
}}}
\put(8431,-2581){\makebox(0,0)[lb]{\smash{\SetFigFont{8}{9.6}{\rmdefault}{\mddefault}{\updefault}{\color[rgb]{0,0,0}$(0,0)$}%
}}}
\put(8041,-2941){\makebox(0,0)[lb]{\smash{\SetFigFont{8}{9.6}{\rmdefault}{\mddefault}{\updefault}{\color[rgb]{0,0,0}$\sigma_2$}%
}}}
\put(8536,-2116){\makebox(0,0)[lb]{\smash{\SetFigFont{8}{9.6}{\rmdefault}{\mddefault}{\updefault}{\color[rgb]{0,0,0}$\sigma_1$}%
}}}
\end{picture}

\end{center}
This gives both the correct affine structure on each two-cell (making
it isomorphic to the standard two-simplex) and the correct
fan structure at the vertices $v_1$ and $v_2$. Note that $\sigma_1$
has the same shape in each chart, so these are identified under the
maps ${\rm exp}_{v_i}$ up to translation. However, up to translation,
the linear transformation $\begin{pmatrix} 1&0\\ 4&1\end{pmatrix}$
is required to transform $\sigma_2$ in $R_{v_1}$ to $\sigma_2$
in $R_{v_2}$. Thus one finds that if one follows the affine coordinates
along a loop starting at $v_2$, into $\sigma_1$ to $v_1$ and
into $\sigma_2$ back to $v_2$, they will undergo an affine transformation
(called the holonomy around the loop) whose linear part is
the linear transformation $\begin{pmatrix} 1&0\\ 4&1\end{pmatrix}$.
Thus in particular, there is no way to extend the affine structure
across the point in $\Delta$ on the line segment joining $v_1$ and $v_2$.

(5) One can carry out this procedure for degenerations of hypersurfaces
or complete intersections in toric varieties. In the hypersurface case,
one obtains the same affine manifolds with singularities described in
\cite{HZ}, or in the case of the quintic, in \cite{GHJ}. 
Details will be given elsewhere.
\end{examples}

\section{From affine manifolds to toric degenerations}

The most important aspect of our proposed construction is the ability
to reverse the construction of the previous section. Given an
integral affine manifold with singularities $B$ and a toric polyhedral
decomposition $\P$, we wish to construct a toric degeneration
$\X\rightarrow S$ coming from this data.

The first step of the construction is easy, i.e. the construction of $\X_0$. 
Let $B,\P$ be as above.
Again, for simplicity, we will assume that no $\sigma\in\P$ is 
self-intersecting:
this is equivalent to the irreducible components
of $\X_0$ being normal. In particular, the endpoints of any
edge in $\P$ are distinct. Let $v_1,\ldots,v_m$ be the vertices
of $\P$. Then the tangent space $\T_{B,v_i}$ contains
a natural integral lattice induced by the integral structure,
(the lattice generated by $\partial/\partial y_1,\ldots,\partial/
\partial y_n$, where $y_1,\ldots,y_n$ are local integral affine
coordinates). Identifying this lattice with $M$ and
$\T_{B,v_i}$ with $M_{\RR}$, and identifying a neighbourhood
of $v_i$ with a neighbourhood of zero in $M_{\RR}$,
$\P$ looks locally near $v_i$ like a fan $\Sigma_i$ in $M_{\RR}$.
This fan defines a toric variety
$X_i$, and $X_1,\ldots,X_m$ will be the irreducible
components of $\X_0$. We then glue together the components of
$\{X_i\}$ using the combinatorics dictated by $\P$. Specifically,
if $v_i$ and $v_j$ are joined by an edge $e\in \P$, then this edge
defines rays in both $\Sigma_i$ and $\Sigma_j$, and hence divisors
$D_i$ and $D_j$ in $X_i$ and $X_j$. It follows from 
condition (5) of Definition \ref{PD2} that $D_i$ and
$D_j$ are isomorphic: they are in fact defined by the same
fan. This isomorphism is canonical, and we glue $X_i$ and
$X_j$ along these divisors using this canonical isomorphism. Again,
one can show that (5) of Definition \ref{PD2} guarantees
all such gluings are compatible, and one obtains a scheme $\X_0$.

There are several points to make here. First, there is actually
a whole moduli space of such gluings. We said above that there
was a canonical isomorphism between $D_i$ and $D_j$, but such
an isomorphism can be twisted by an automorphism of $D_i$. Thus
if one specifies, for each pair $i,j$ an automorphism of $D_i$,
and demands in addition some compatibility conditions, one obtains
a new gluing which is not necessarily isomorphic to the original
gluing. In fact, one can parametrize the set of all possible gluings
by a \v Cech cohomology group of a sheaf on $B$.

The second point is that just knowing $\X_0$ gives nowhere near
enough information to smooth $\X_0$ correctly. In particular, we
have so far only used the fan structure of $\P$,
and used no information about the maximal cells of $\P$. 
There may be different smoothings of $\X_0$ depending on this data:
we saw this in Example \ref{Bexamp}, (2) and (3), where there can be many
different smoothings of a cycle of rational curves, giving different
singularities of the total space of the smoothing. In particular, 
just knowing $\X_0$ tells us nothing about the toric varieties
$Y_x$ which may appear as local models for a smoothing
$\X\rightarrow S$. We rectify
this by introducing {\it log structures}.

Roughly put, a log structure is some additional structure on $\X_0$
which reflects some essentially toric information about the embedding
$\X_0\subseteq\X$. This may be viewed as something akin to
infinitesimal information about the smoothing.

Recall a {\it monoid} is a set with an associative product with
a unit. We will only use commutative monoids here.

\begin{definition}
A log structure on a scheme (or analytic space) $X$ is a (unital) homomorphism
$$\alpha_X:\shM_X\rightarrow \O_X$$
of sheaves of (multiplicative) monoids inducing an isomorphism
$\alpha_X^{-1}(\O_X^{\times})\rightarrow \O_X^{\times}$. The
triple $(X,\shM_X,\alpha_X)$ is then called a {\it log space}.
We often write the whole package as $X^{\dagger}$.
\end{definition}

A morphism of log spaces $F:X^{\dagger}\rightarrow Y^{\dagger}$ consists
of a morphism $\underline{F}:X\rightarrow Y$ of underlying
spaces together with a homomorphism $F^{\#}:\underline{F}^{-1}(\shM_Y)
\rightarrow\shM_X$ commuting with the structure homomorphisms:
$$\alpha_X\circ F^{\#}=\underline{F}^*\circ\alpha_Y.$$

The key examples:

\begin{examples}
\label{logexamples}
(1) Let $X$ be a scheme and $D\subseteq X$ a closed subset of
codimension one. Denote by $j:X\setminus D\rightarrow X$
the inclusion. Then the inclusion
$$\alpha_X:\shM_X=j_*(\O_{X\setminus D}^{\times})\cap\O_X\rightarrow
\O_X$$
of the sheaf of regular functions with zeroes contained in $D$ is a log
structure on $X$.

(2) A {\it prelog structure}, i.e. an arbitrary homomorphism of
sheaves of monoids $\varphi:\shP\rightarrow\O_X$, defines
an associated log structure $\shM_X$ by
$$\shM_X=(\shP\oplus\O_X^{\times})/\{(p,\varphi(p)^{-1})|p\in
\varphi^{-1}(\O_X^{\times})\}$$
and $\alpha_X(p,h)=h\cdot\varphi(p)$.

(3) If $f:X\rightarrow Y$ is a morphism of schemes and $\alpha_Y:\shM_Y
\rightarrow\O_Y$ is a log structure on $Y$, then the prelog structure
$f^{-1}(\shM_Y)\rightarrow\O_X$ defines an associated log structure
on $X$, the {\it pull-back log structure}.

(4) In (1) we can pull back the log structure on $X$ to $D$ using
(3). Thus in particular, if $\X\rightarrow S$ is a toric
degeneration, the inclusion $\X_0\subseteq\X$ gives a log
structure on $\X$ and an induced log structure on $\X_0$. Similarly
the inclusion $0\in S$ gives a log structure on $S$ and
an induced one on $0$. Here $\M_0=\CC^{\times}\oplus\NN$,
where $\NN$ is the (additive) monoid of natural (non-negative) numbers,
and 
$$\alpha_0(h,n)=\begin{cases}h& n=0\\ 0&n\not=0.\end{cases}$$
We then have log morphisms $\X^{\dagger}\rightarrow S^{\dagger}$ and
$\X_0^{\dagger}\rightarrow 0^{\dagger}$.

(5) If $\sigma\subseteq M_{\RR}$ is a cone, $\dual{\sigma}\subseteq
N_{\RR}$ the dual cone, let $P=\dual{\sigma}\cap N$: this is a monoid.
The affine toric variety defined by $\sigma$ can be written as 
$X=\Spec \CC[P]$. Here $\CC[P]$ denotes the monoid ring of $P$,
generated as a vector space over $\CC$ by symbols $\{z^p|p\in P\}$
with multiplication given by $z^p\cdot z^{p'}=z^{p+p'}$.

We then have a pre-log structure induced by the homomorphism of
monoids 
$$P\rightarrow \CC[P]$$
given by $p\mapsto z^p$. There is then an associated log
structure on $X$. If $p\in P$, then the monomial $z^p$ defines a map
$f:X\rightarrow \Spec \CC[\NN]\quad (=\Spec \CC[t])$ which is a log morphism.
The fibre $X_0=\Spec \CC[P]/(z^p)$ is a subscheme of $X$,
and there is an induced log structure on $X_0$, and a map $X_0^{\dagger}
\rightarrow 0^{\dagger}$ as in (4).

Condition (4) of Definition \ref{toricdegen} in fact implies
that locally, away from $Z$, $\X^{\dagger}$ and $\X_0^{\dagger}$ are
of the above form.
\end{examples}

\begin{remark}
It is sometimes useful to think about a log structure via the exact
sequence
\begin{eqnarray}
\label{exactseq}
1\rightarrow\O_X^{\times}\rightarrow\shM_X\rightarrow\overline{\shM}_X
\rightarrow 0
\end{eqnarray}
defining a sheaf of monoids $\overline{\shM}_X$. For example, consider
Example \ref{logexamples}, (1), with 
$D=\{x_1x_2=0\}\subseteq X=\Spec k[x_1,x_2]$, $D=D_1\cup D_2$,
with $D_i=\{x_i=0\}$. If $i_j:D_j
\rightarrow X$ are the inclusions, then $\overline{\shM}_X
=i_{1*}\NN\oplus i_{2*}\NN$, and an element $f\in\shM_X$ is mapped
to $(n_1,n_2)$, where $n_j$ is the order of vanishing of $f$
along $D_j$. Pulling back this log structure to $D$, one
obtains a similar exact sequence with $\overline{\shM}_D=\overline{\shM}_X$.
\end{remark}

Example \ref{logexamples}, (5) is important. 
The beauty of log geometry is that 
we are able to treat such an $\X_0^{\dagger}$ as if it were
a non-singular variety. Essentially, we say a log scheme over $0^{\dagger}$
which is
locally of the form given in (5) is {\it log smooth}
over $0^{\dagger}$. In particular, on a log smooth scheme,
there is a sheaf of logarithmic differentials which is locally
free. F. Kato has developed deformation theory for such log schemes.

The philosophy is then as follows. Given $B,\P$, we have constructed
a space $\X_0$. We first try to put a log structure on $\X_0$
such that there is a set $Z\subseteq Sing(\X_0)$ not
contained in any toric stratum of $\X_0$ such that $\X_0^{\dagger}
\setminus Z$ is log smooth. We then try to deform $\X_0^{\dagger}$
in a family to get $\X^{\dagger}\rightarrow S^{\dagger}$,
and do this in such a way that the underlying map of spaces
$\X\rightarrow S$ is a toric degeneration.

This is the technical heart of the program, and we will only give some
hints here of how this works. 

The first point is that given $B,\P$, we constructed $\X_0$
by gluing together its irreducible components. However, it can
also be constructed by describing an open cover and a gluing of these
open sets. Specifically, if $\sigma\in\P$, we can view $\sigma$
as a polytope in $M_{\RR}$, and let $\tilde\sigma$ be the cone
over $\sigma$ in $M_{\RR}\oplus\RR$, i.e.
$$\tilde\sigma=\{(rm,r)|r\in \RR_{\ge 0}, m\in\sigma\}.$$
Then $\tilde\sigma$ defines an affine toric variety $Y_{\sigma}$,
and writing $\Hom_{\ZZ}(M\oplus\ZZ,\ZZ)$ as $N\oplus\ZZ$, $\rho=(0,1)
\in N\oplus \ZZ$ represents a monomial $z^{\rho}$ on $Y_{\sigma}$.
We let $X_{\sigma}$ be defined by the equation $z^{\rho}=0$ in $Y_{\sigma}$.
We are simply reversing the procedure described in \S 2 to obtain $B$
from a toric degeneration. It is then not difficult to show
that $\{X_{\sigma}|\sigma\in\P\}$ form a natural open covering of
$\X_0$, and one can explicitly describe how they glue.

We have gained from this description a natural log structure on $X_{\sigma}$
coming from the inclusion $X_{\sigma}\subseteq Y_{\sigma}$ as
in Example \ref{logexamples}, (4) or (5), and so we have an
open covering of $\X_0$ by log schemes which are log smooth
over $0^{\dagger}$, $\{X_{\sigma}^{\dagger}\rightarrow 0^{\dagger}|
\sigma\in\P\}$. 

The problem is that unless $Z=\phi$, the log structures don't glue.
However, it is possible to define a ``sheaf of log structures on $X_{\sigma}$
over $0^{\dagger}$'' consisting essentially of deformations of the
given log structure. This is a rather technical but important point.

Given a log structure, 
in a certain sense,
it is the extension class of $\overline{\shM}_X$ by $\O_X^{\times}$
in (\ref{exactseq})
which determines the log structure. In our case, the sheaves
$\overline{\shM}_{X_{\sigma}}$ do glue, so one can define a sheaf
of monoids $\overline{\shM}_{\X_0}$ globally on $\X_0$. The ``sheaf of
log smooth structures'' $LS(\X_0)$ can then be defined as an appropriate
subsheaf of $\shExt^1(\overline{\shM}_{\X_0}^{gp},\O_{\X_0}^{\times})$.
Here the superscript $gp$ refers to the Grothendieck group of the monoid.

\begin{examples}
(1) A simple example shows that one can have a non-trivial family
of log structures even very locally: let $X\subseteq \CC^4=\Spec
\CC[x,y,z,w]$ be given by the equations $xy=zw=0$. Then there is a natural
one-parameter family of log structures on $X$ induced by the inclusions
$X\subseteq Y_{\lambda}$, where $Y_{\lambda}$ is the quadric given
by the equation $xy-\lambda zw=0$, $\lambda\in\CC^{\times}$. This gives a family of log structures
$X_{\lambda}^{\dagger}$, and there is no isomorphism $f:X^{\dagger}_{\lambda}
\rightarrow X^{\dagger}_{\lambda'}$ of log schemes which is the
identity on the underlying scheme $X$ unless $\lambda=\lambda'$.

(2) If $\X_0^{\dagger}$ 
has normal crossings, (i.e. all elements of $\P$ are
standard simplices) and if $D=Sing(\X_0)$, there is a standard line
bundle $\shN_D$ on $D$ which can be defined as the sheaf of local
infinitesimal deformations $\shExt^1(\Omega^1_{\X_0},\O_{\X_0})$.
Then $LS(\X_0)$ turns out to be the $\O_D^{\times}$-torseur associated
to $\shN_D$, i.e. the sheaf of nowhere zero sections of $\shN_D$.
Thus there only exists a log smooth structure on $\X_0$ if $\shN_D\cong
\O_D$ (see \cite{KN}). In fact, this will only be the case if the minimal
discriminant locus $\Delta\subseteq B$ is empty, and one can
read off $\shN_D$ from information about the discriminant locus.
This is the reason we must allow singularities in the log structure
in the presence of singularities on $B$.

(3) The reader may wonder why we can't restrict to normal crossings
outside of $Z$, where from (2) it appears the theory is relatively
simple.

There are several reasons why we can't do this and don't want to do this.
Indeed, the toric 
situation is the most natural one. First, one might argue you
could try to further subdivide $\P$ so it consists only of standard
simplices. Even if $\Delta$ is empty, this cannot be done
in dimension four or higher, as it is well-known that in these
dimensions there exists simplices which cannot be subdivided
into elementary ones. More seriously, one can't even always subdivide
into simplices. In three dimensions, say, if one has a long strand
in $\Delta$, one finds the definition of a toric polyhedral decomposition
requires that any one-dimensional cell intersecting this strand
is constrained to lie in
a given plane containing that strand and parallel to a certain line
contained in that plane. This results in decompositions which involve
polyhedra which necessarily aren't simplices, as depicted in the following
picture:
$$\epsfbox{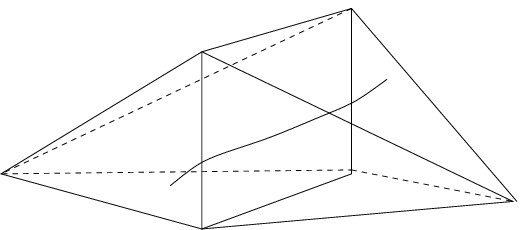}$$

Finally, even if we do have a normal crossings example, the mirror
construction given in the next section will almost certainly never
result in a mirror family which is normal crossings. So the toric
setting is imperative for this approach.
\end{examples}

Thus, in general, we can obtain log smooth structures on dense open subsets
$U\subseteq\X_0$, with $\X_0\setminus U=Z$ the singular set.
However, there are additional restrictions on the nature of 
the singularities of the log structure to ensure there is a local
smoothing. For example, in the normal crossings case, we can
only allow sections of $LS(\X_0)$ to have zeroes, but not poles,
along the singular set $Z$; i.e. the allowable log structures on
$\X_0$ are determined by sections of $\shN_D$. Thus some additional
positivity condition is required, and this is reflected in the
holonomy of the affine structure about $\Delta$.

This allows us to construct a moduli space of allowable log structures
on $\X_0$. Combining this with the fact that $\X_0$ itself might
have locally trivial deformations coming from a family of possible
gluings of the irreducible components, we obtain a whole moduli space
of such log schemes.

It remains an open question being actively researched to determine
when such $\X_0^{\dagger}$ is smoothable. This will be the final
step of the program. We restrict our comments here to the statement
that in dimension $\le 2$, suitably positive $\X_0^{\dagger}$ are
always smoothable, but in dimension $\ge 3$, one can construct
non-smoothable examples (analagous to examples of non-smoothable
three-dimensional canonical singularities). However, we expect
smoothability to be implied by certain conditions on $B$.

\section{The intersection complex, the discrete Legendre transform, and
mirror symmetry}

Now let $\X\rightarrow S$ be a toric degeneration, $\dim\X_0=n$. Suppose in
addition it is polarised by a relatively ample line bundle, i.e. a line
bundle $\shL$ on $\X$ which restricts to an ample line bundle on $\X_t$
for each $t\in S$ (including $t=0$). Then one can construct
another integral affine manifold with singularities $\check B$ with
a polyhedral decomposition $\check\P$ from this data. Just as $B,\P$
of \S 3 was the dual intersection complex, $\check B,\check\P$ is
the intersection complex, but the polarization is necessary to define
the affine structure.

There are two ways to think about $\check B$. First, the more direct
way is as follows. Again assume every component of $\X_0$
is normal. For each component $X_i$ of $\X_0$, let
$\sigma_i$ be the Newton polytope of $\shL|_{X_i}$ in $N_{\RR}$. This
is a lattice polytope which is well-defined up to affine linear 
transformations. The cells of $\sigma_i$ are in a one-to-one
inclusion preserving correspondence with the toric strata of $X_i$,
and we can then identify cells of $\sigma_i$ and $\sigma_j$ if
the corresponding toric strata of $X_i$ and $X_j$ are identified
in $\X_0$. Making these identifications gives $\check B$.
To get the singular affine structure on $\check B$, we need to specify
a fan structure at each vertex of $\check B$. But a vertex of $\check B$
corresponds
to $\{x\}\in Strata(\X_0)$ with an associated polyhedron $\sigma_x\subseteq
M_{\RR}$. We define the normal fan $\check\Sigma_x$ of $\sigma_x$
to be the fan whose cones are in
one-to-one inclusion reversing
correspondence with the cells $\tau$ of $\sigma_x$, with the cone
corresponding to $\tau$ given by
$$\{f\in N_{\RR}|
\hbox{$f|_{\tau}$ is constant and $\langle f,z\rangle\ge
\langle f,y\rangle$ for all $z\in\sigma_x,y\in\tau$}\}$$
We then take the fan $\check\Sigma_x$ to determine the fan structure
at the corresponding vertex of $\check B$.

\begin{proposition}
\label{intersectioncomplex}
$\check B$ is an integral affine manifold with singularities, and $\check
\P$ is a toric polyhedral decomposition on $\check B$. If 
$\check\Delta\subseteq\check B$ is the minimal discriminant locus
of the affine structure, then there is a homeomorphism $\alpha:B
\rightarrow\check B$ with $\alpha(\Delta)=\check\Delta$. Furthermore,
the affine structures on $B$ and $\check B$ are dual in the sense
that the holonomy representations of the flat connections on $\T_B$
and $\T_{\check B}$ induced by the respective affine structures are
naturally dual.
\end{proposition}

$\check B$ can be described more intrinsically in terms of $B$ and
$\P$ using the discrete Legendre transform.

\begin{definition}
\label{PLfunc} If $B$ is an affine manifold with singularities,
let $\Aff(B,\RR)$ denote the sheaf of functions on $B$ with values
in $\RR$ which are affine linear when restricted to $B_0$.
Let $\P$ be a polyhedral decomposition of $B$. If $U\subseteq B$ is an open set,
then a {\it piecewise linear} function on $U$
is a continuous function $f:U\rightarrow \RR$ which is affine linear on $U
\cap Int(\sigma)$ for each maximal $\sigma\in\P$, and which satisfies
the following property: for any $y\in U$, $y\in Int(\sigma)$ for some
$\sigma\in\P$, there exists a neighbourhood $V$ of $y$ and
a $g\in \Gamma(V,\Aff(B,\RR))$ such that $f-g$ is zero on
$V\cap Int(\sigma)$.
\end{definition}

\begin{definition}
A {\it multi-valued piecewise linear function} on $B$
with respect to $\P$ is a collection of piecewise linear functions
$(U_i,\varphi_i)$ for $\{U_i\}$ an open cover of $B$, such that
$\varphi_i-\varphi_j\in \Gamma(U_i\cap U_j,\Aff(B,\RR))$ for all $i,j$.
We say a multivalued piecewise linear function $\varphi$ is {\it
strictly convex} if at every vertex $v$ of $\P$, some representative
$\varphi_i$ is strictly convex in a neighbourhood of $v$ (in the usual
sense of a strictly convex piecewise linear function
on a fan).
\end{definition}

 Let $\varphi$ be a strictly convex multi-valued piecewise linear
function on $B$ with only integral slopes. We will construct a new integral
affine manifold with singularities $\check B$ with discriminant locus
$\check\Delta$. As manifolds, $B=\check B$ and $\Delta=\check\Delta$,
but the affine structures are dual. In addition, we obtain a toric polyhedral
decomposition $\check\P$ and a strictly convex multi-valued piecewise
linear function $\check\varphi$ with integral slope on $\check B$. We will
say $(\check B,\check\P,\check\varphi)$ is the {\it discrete Legendre transform}
of $(B,\P,\varphi)$.

First we define $\check\P$. For any $\sigma\in\P$, define $\check\sigma$
to be the union of all simplices in $Bar(\P)$ intersecting $\sigma$ but
disjoint from any proper subcell of $\sigma$. Put
$$\check\P=\{\check\sigma|\sigma\in\P\}.$$
This is the usual dual cell complex to $\sigma$, with $\dim\check\sigma=n-
\dim\sigma$.
Of course, $\check\sigma$ is not a polyhedron with respect to the affine
structure on $B$, and we will build a new affine structure on $B$ using 
the method of \S 2.

For any vertex $v\in\P$, we obtain a fan $\Sigma_v$ living in $\T_{B,v}$,
and locally, $\varphi$ defines a piecewise linear function $\varphi_v$
on the fan $\Sigma_v$ up to a choice of a linear function.
This function is strictly convex by assumption, and we can
consider the corresponding Newton polytope, i.e. set
$$\check v'=
\{x\in\T_{B,v}^*|\langle x,y\rangle
\ge -\varphi_v(y)\quad \forall y\in\T_{B,v}\}.$$
Note that because $\varphi_{v}$ is strictly convex there is a
one-to-one inclusion
reversing correspondence between the cells of $\check v'$
and cones in $\Sigma_{v}$; if $\tau\in\Sigma_{v}$, the corresponding
cell $\check\tau\subseteq\check v'$ is
$$\check\tau=\{x\in\check v'|\langle x,y\rangle=-\varphi_{v}(y)\quad
\forall y\in\tau\}.$$
In addition, $\check v'$ is an integral polytope because $\varphi_{v}$
has integral slopes.

Each $\check v'$ can then be identified in a canonical way
with $\check v\in\check \P$.
This can be done in a piecewise linear way on each simplex of
the first barycentric subdivision of $\check v'$. This gives an identification
of $\check v$ with a lattice polytope in $N_{\RR}$, giving the first
step of the construction of the dual affine structure on $B$.

To finish specifying an integral affine structure with singularities
on $\check B=B$,
we just need to specify a fan structure at
each vertex $\check\sigma$ of $\check\P$ (for $\sigma$ a maximal cell
of $\P$). We take the fan structure at $\check\sigma$ to be
given by the normal fan 
$\check\Sigma_{\sigma}$ of $\sigma$, just as before.

Finally, we wish to define $\check\varphi$, the Legendre transform
of $\varphi$. We do this by defining $\check\varphi$ in a neighbourhood
of each vertex $\check\sigma$ of $\check\P$, where $\sigma\in\P$
is a maximal cell. This is equivalent to giving a piecewise linear
function $\check\varphi_{\check\sigma}$ on the normal fan 
$\check\Sigma_{\sigma}$ of $\sigma$, viewing $\sigma$ as a polytope in
$M_{\RR}$. Since we want the operation of discrete Legendre transform
to be a duality, $\sigma$ must be obtained as the Newton polytope
coming from the function $\check\varphi_{\check\sigma}$ on 
$\check\Sigma_{\sigma}$, and thus we are forced to define 
$\check\varphi_{\check\sigma}$ by
$$\check\varphi_{\check\sigma}(y)
=-\inf\{\langle y,x\rangle|x\in\tilde\sigma\}$$
for $y\in N_{\RR}$. This is a piecewise linear function
on the fan $\check\Sigma_{\sigma}$, and it is a
standard easy fact that it is strictly convex, with the Newton polyhedron
of $\check\varphi_{\check\sigma}$ being $\sigma$. If $\sigma$
is shifted in $M_{\RR}$ by a translation, 
$\check\varphi_{\check\sigma}$ is changed
by a linear function, so it is well-defined modulo linear functions.

Thus given the triple $(B,\P,\varphi)$, we obtain $(\check B,\check\P,
\check\varphi)$. 

Now the point is that given a toric degeneration $\X\rightarrow S$,
polarized by $\shL$, we actually obtain a strictly convex
multi-valued piecewise linear function on $B$: in a neighbourhood
of each vertex $v_i$ of $\P$ corresponding to 
the irreducible component $X_i$ of $\X_0$, the line bundle $\shL|_{X_i}$
yields a piecewise linear function, up to a linear function, on
the corresponding fan $\Sigma_i$. This defines a piecewise linear function
$\varphi_i$ in a neighbourhood of $v_i$ on $B$. One
can check these define a multi-valued piecewise linear function on $B$,
and it is strictly convex because $\shL|_{X_i}$ is ample for each $i$.
Thus we obtain a triple $(B,\P,\varphi)$, which we call the {\it degeneration
data} associated to $\X\rightarrow S$, $\shL$. 

It is easy to see that the first construction of $\check B,\check\P$
given here as the intersection complex
coincides with the data from $(\check B,\check\P,\check\varphi)$,
the discrete Legendre transform of $(B,\P,\varphi)$.

We now come to the fundamental idea of the paper: {\it mirror
symmetry can be understood as a duality between toric degenerations
with dual degeneration data}. For example, if $f:\X\rightarrow S$
and $\check f:\check\X\rightarrow S$ are polarized toric degenerations
such that their degeneration data $(B,\P,\varphi)$ and $(\check B,\check\P,
\check\varphi)$ are related by the discrete Legendre transform, then $f$
and $\check f$ should be viewed as mirror degenerations. One should also 
consider the singular fibres by themselves: the log schemes $\X_0^{\dagger}$
and $\check\X_0^{\dagger}$ along with polarisations carry enough data
by themselves to define $(B,\P,\varphi)$ and $(\check B,\check\P,\check
\varphi)$. We can then say $\X_0^{\dagger}$ and $\check\X_0^{\dagger}$ are
a mirror pair of log schemes if again the degeneration data are related by the 
discrete Legendre transform.

Strictly speaking, mirror symmetry should be about families. One can make
more precise statements, defining suitable moduli spaces of log schemes and 
log K\"ahler moduli. Under suitable hypotheses, roughly implying
the degeneration is a large complex structure limit
rather than just a maximally unipotent degeneration, one can define
a natural mirror map identifying these two different moduli spaces for
a mirror pair. Furthermore, in this case, we expect that we can
deduce results about the smoothings $\X_t$ and $\check\X_t$ for
$t\not=0$, and in particular demonstrate that 
$h^{1,1}(\X_t)=h^{1,n-1}
(\check\X_t)$ and $h^{1,n-1}(\X_t)=h^{1,1}(\check\X_t)$. It is likely,
given appropriate assumptions, that this will follow from standard
techniques of log geometry. Of course, it cannot hold in general because
$\X_t$ need not be smooth, but only have canonical singularities, in which
case it is not clear what the correct equalities should be.

This conception of mirror symmetry fits with the Batyrev-Borisov mirror symmetry
construction \cite{BatBor}, and generalises that construction. While it
also applies to degenerations of complex tori, where traditional
forms of mirror symmetry for complex tori are reproduced, it is
not clear how much broader this construction is. However, we expect
it should be significantly broader, and it certainly puts many
different forms of mirror symmetry on an equal footing.

Philosophically, once the details of the basic construction are complete,
one can hope that one can study mirror symmetry for non-singular Calabi-Yau
manifolds by studying mirror symmetry for log schemes. Many objects of interest
should have analogous log versions. For example, if one is interested in 
computing Gromov-Witten invariants, one can try to define log Gromov-Witten
invariants, which can be computed on the singular fibre $\X_0^{\dagger}$.
If defined appropriately, these invariants will remain stable under smoothing,
and so one reduces the calculation of Gromov-Witten invariants to the
singular case, which may be easier. In particular, one should be able
to relate such calculations to questions of combinatorics of graphs on
$B$, much as in \cite{KS} and \cite{Fukaya}. 
This approach has been started by the second
author in \cite{SlGW}, and though technical problems remain, we believe
log Gromov-Witten invariants can be defined. Previous general
work in this direction is due to Tian; Li and Ruan \cite{LR};
Ionel and Parker \cite{IP}; Gathmann \cite{Gath} and Li \cite{Li1},\cite{Li2}, 
who
covers cases relative a smooth divisor and of normal crossing varieties with
smooth singular locus.
Nevertheless, this approach remains a major undertaking.

\begin{example} Let $f:\X\rightarrow S$ be a degeneration
of elliptic curves as in Example \ref{Bexamp}, (3), with $B=\RR/m\ZZ$,
decomposed into line segments of lengths $m_1,\ldots,m_p$ with $\sum m_i=m$,
so that $\X_0$ has $p$ components. Choose a polarization on $\X$ of
degree $n$, which is degree $n_1,\ldots,n_p$ on the $p$ components
of $\X_0$ respectively ($n_i\ge 1$, $\sum n_i=n$). Then $\check B$
is a union of line segments of lengths $n_1,\ldots,n_p$, so $\check B
=\RR/n\ZZ$. This yields a new degeneration $\check f:\check\X
\rightarrow S$, with a polarization which is degree $m_1,\ldots,m_p$
on the irreducible components of $\check\X_0$.

Of course, the precise geometry of the singular fibres and their
polarizations are dependent on the initial choice of polyhedral
decomposition $\P$ and $\varphi$, but if one deletes the singular
fibres, these differences disappear. Letting $S^*=S\setminus\{0\}$,
$f:\X^*=f^{-1}(S^*)\rightarrow S^*$, and $\check f:\check\X^*\rightarrow S^*$
give two families of polarized elliptic curves, one with monodromy in
$H^1$ being $\begin{pmatrix} 1&m\\ 0&1\end{pmatrix}$ in a suitable basis
of cohomology, with the polarization of degree $n$, and the other 
with the roles of $m$ and $n$ reversed. This is the manifestation
of mirror symmetry here.
\end{example}

\begin{remark} This construction should give mirror pairs for toric
degenerations, at least with certain additional conditions, but
there remains the question of how general a mirror symmetry
construction this is. More generally, the general mirror symmetry
conjecture suggests there should be mirror partners associated to any
large complex structure limit (see \cite{Morlcsl}). Now in general
a toric degeneration is a maximally unipotent degeneration, but
does not necessarily satisfy the stronger condition of being
a large complex structure limit. We do not expect that any
maximally unipotent degeneration is birationally equivalent to a toric
degeneration. However, there is some tenuous evidence which leads
one to speculate that any large complex structure degeneration is in fact
birationally equivalent to a toric degeneration. If this is the case,
our proposed construction will yield a general mirror symmetry
construction.
\end{remark}

\section{Connections with the Strominger-Yau-Zaslow conjecture}

This approach to mirror symmetry is closely related to the 
Strominger-Yau-Zaslow approach, and we believe that it should be viewed
as an algebro-geometrisation of SYZ, a discretization, so to speak.

Recall briefly that in Hitchin's approach to SYZ \cite{Hit}
(see also \cite{Leung}), 
one considers an affine manifold $B$ whose transition functions are
contained in $\Aff_{\RR}(M)$ (rather than $\Aff(M)$ or $\Aff(M_{\RR})$).
Then there is a well-defined family of sublattices $\Lambda$ of $\T_B$
generated by $\partial/\partial y_1,\ldots,\partial/\partial y_n$
if $y_1,\ldots,y_n$ are local affine coordinates. Because of the 
restriction on transition functions, this basis is well-defined
up to elements of $GL_n(\ZZ)$; hence $\Lambda\subseteq\T_B$ is well-defined.
We can set $X(B)=\T_B/\Lambda$. This comes along with a complex
structure: locally one can define holomorphic coordinate functions
$z_i$ on $X(B)$ over a point in $B$ with coordinates $(y_1,\ldots,y_n)$,
$$z_i\left(\sum_j x_j{\partial\over\partial y_j}\right)=
e^{2\pi\sqrt{-1}(x_j+\sqrt{-1}y_j)}.$$
If we also have a strictly convex differentiable function
$\varphi$ on $B$ (i.e. $(\partial^2\varphi/\partial y_i\partial y_j)$ is
positive definite for affine coordinates $y_1,\ldots,y_n$),
then if $\pi:X(B)\rightarrow B$ is the projection, 
$\varphi\circ\pi$ is the K\"ahler potential of a K\"ahler metric
on $X(B)$ (Ricci flat if $\varphi$ satisfies the real Monge-Amp\`ere
equation $\det(\partial^2\varphi/\partial y_i\partial y_j)=constant$).

To obtain the mirror of $X(B)$, one defines a new affine structure
on $B$ with local coordinates given by $\check y_i=\partial\varphi/\partial
y_i$. One also obtains a function
$$\check\varphi(\check y_1,\ldots,\check y_n)
=\sum_{i=1}^n y_i\check y_i-\varphi(y_1,\ldots,y_n),$$
the Legendre transform of $\varphi$.

The data $\check B,\check\varphi$ defines a new K\"ahler manifold $X(\check B)$
which is SYZ-dual to $X(B)$.

One of the difficulties with this approach is that this is only an
approximation, and it only gives mirror symmetry precisely in
the complex torus case. In other cases, one expects singular fibres which
destroy the ability to have a pleasant complex structure on $X(B)$
(and in fact if $B$ has singularities, we don't know how to define
$X(B)$). However, one can show these bundles can
be good approximations to the genuine complex structures.

Let us move towards a more precise statement here. First, we need to
underline the importance of {\it integral} affine structures for
the study of large complex structure limits. This was first
observed by Kontsevich and Soibelman \cite{KS}. For convenience, let us
first introduce a few additional concepts of affine manifolds and torus
bundles over them.

Let $\pi:\tilde B\rightarrow B$ be the universal
covering of an (integral)
affine manifold $B$, inducing an (integral) affine structure on
$\tilde B$. Then there is an (integral) 
affine immersion $d:\tilde B\rightarrow M_{\RR}$,
called the {\it developing map}, and any two such maps differ only
by an (integral) affine transformation. This map is obtained in a
standard way by patching together (integral) affine coordinate charts
on $B$. We can then obtain the {\it holonomy representation} as follows.
The fundamental group $\pi_1(B)$ acts on $\tilde B$ by deck transformations; for
$\gamma\in\pi_1(B)$, let $T_{\gamma}:\tilde B\rightarrow\tilde B$
be the corresponding deck transformation with $T_{\gamma_1}\circ
T_{\gamma_2}=T_{\gamma_2\gamma_1}$.
Then by the uniqueness of
the developing map, there exists a $\rho(\gamma)\in \Aff(M_{\RR})$
such that $\rho(\gamma)\circ d\circ T_{\gamma}=d$. The map
$\rho:\pi_1(B)\rightarrow \Aff(M_{\RR})$ is called the {\it holonomy
representation}. If the affine structure is integral, then $\im\rho
\subseteq \Aff(M)$.

Now suppose $B$ is equipped with an integral affine structure,
with developing map $d:\tilde B\rightarrow M_{\RR}$ and holonomy
representation $\rho:\pi_1(B)\rightarrow \Aff(M)$. Suppose $d':\tilde
B\rightarrow M_{\RR}$ is another map, not necessarily an immersion, such
that there is a representation $\rho':\tilde B\rightarrow \Aff(M_{\RR})$
satisfying $\rho'(\gamma)\circ d'\circ T_{\gamma}=d'$. Assume that $\rho$
and $\rho'$ have the same linear parts, i.e. agree if composed with the
natural projection $\Aff(M_{\RR})\rightarrow GL(M_{\RR})$. 
Assume furthermore that for $r>R$, $rd+d':\tilde B\rightarrow M_{\RR}$ is
an immersion. Then in
this case we can form a new affine manifold
$$\bar B=B\times (R,\infty)$$
defined by a developing map
$$\bar d:\tilde B\times (R,\infty)\rightarrow M_{\RR}\times\RR$$
given by
$$\bar d(b,r)=(rd(b)+d'(b),r).$$
One can check that this defines an affine structure on $\bar B$ with
transition maps (or holonomy representation) contained in 
$$(M_{\RR}\oplus\RR)\rtimes GL(M\oplus\ZZ),$$
i.e. has integral linear part. If the affine structure coming from
$d$ had not been integral to begin with, we would not have had
the integrality of the linear part for $\bar B$, and thus been unable to
form the complex manifold $X(\bar B)$.

$X(\bar B)$ is a complex manifold of dimension $n+1$. The projection
$\bar B\rightarrow (R,\infty)$ induces a map of complex manifolds
$X(\bar B)\rightarrow S^*=X((R,\infty))$, the latter being a punctured
disk of radius $e^{-2\pi R}$. This is a family of $n$-dimensional
complex manifolds, with, in general, non-trivial monodromy.

There is one more refinement of this. Given an open covering $\{U_i\}$
of $B$ and $\alpha=(\alpha_{ij})$ a \v Cech 1-cocycle of flat
sections of $\T_B/\Lambda$, we can glue $X(U_i)$ and $X(U_j)$ by fibrewise
translating
by $\alpha_{ij}$ before gluing. This translation is holomorphic, and so
the glued manifold is still a complex manifold, which can be viewed
as a twisted form of $X(B)$. We write this as $X(B,\alpha)$.

Next we define what we mean by a small deformation.

For an affine manifold $B$ with the linear part of the holonomy
being integral, we would like to define the
notion of a small deformation of the complex manifold $X(B,\alpha)$.
First, we need to answer the question of how we should think of a 
deformation of this manifold. One begins with a fixed covering
$\{U_i\}$ of $B$ along with affine coordinate charts
$\psi_i:U_i\rightarrow M_{\RR}$ which are open immersions. Thus there is 
a natural identification of $X(U_i)$ with an open subset
of $X(M_{\RR})\cong M_{\RR}\otimes\CC^{\times}$. Furthermore
$X(B,\alpha)$ is obtained by gluing together the sets
$X(U_i)$ and $X(U_j)$ (for each $i$ and $j$) along
$X(U_i\cap U_j)$ using a biholomorphic map
$\varphi_{ij}:X(U_i\cap U_j)\rightarrow X(U_i\cap U_j)$.
(These maps depend on $\alpha$; if $\alpha=0$ they are the identity.)

To deform $X(B,\alpha)$, we should simply perturb these
maps. To measure the size of this deformation, we measure the size of the
perturbation. First, to measure distance between two points $x,y$
in $X(M_{\RR})=M_{\RR}+i M_{\RR}/M$, take
$$d(x,y)=\inf \{\|\tilde x-\tilde y\| |\hbox
{$\tilde x$, $\tilde y$ are lifts of $x,y$ to $M_{\RR}+i M_{\RR}$}\}.$$
Here $\| \cdot \|$ denotes the norm with respect to some fixed inner
product on $M_{\RR}\otimes\CC$. 

Formally, we say a complex manifold $X$ is a {\it deformation of $X(B,
\alpha)$ of size $C$} 
if $X$ can be covered by open sets $X_i$ along with isomorphisms
$\varphi_i:X(U_i)\rightarrow X_i$ such that
\begin{itemize}
\item
For any point $x\in\varphi_i^{-1}(X_i\cap X_j)$, there is a point
$y\in X(U_i\cap U_j)$ such that $d(x,y)<C$, and
conversely for any $y\in X(U_i\cap U_j)$, there exists an 
$x\in\varphi_i^{-1}(X_i\cap X_j)$ such that
$d(x,y)<C$. (In other words, the gluing sets have only changed by distance
at most $C$).
\item
For any point
$x_i\in\varphi_i^{-1}(X_i\cap X_j)$ and
$x_j\in X(U_j)$ with $\varphi_j(x_j)=\varphi_i(x_i)$,
there exists a point $y\in X(U_i\cap U_j)$ with $d(x_i,y)<C$ and
$d(x_j,\varphi_{ij}(y))<C$.
\end{itemize}

One can then show the following.

\begin{theorem} Let $f:\X\rightarrow S$ be a toric degeneration,
$B$ the corresponding integral affine manifold with singularities.
Then there exists
\begin{itemize} 
\item an open set $U\subseteq B$ such that $B\setminus U$ retracts onto
the singular set $\Delta$.
\item A map $d':\tilde U\rightarrow M_{\RR}$ as above defining
an affine structure on $\bar U=U\times (R,\infty)$.
\item A \v Cech 1-cocycle $\alpha$ of flat sections of $\T_{\bar U}$,
hence giving a map
$$g:X(\bar U,\alpha)\rightarrow X((R,\infty))=\{t\in\CC| 0<|t|
<e^{-2\pi R}\}.$$
\item An identification of $X((R,\infty))$ with a punctured open
neighbourhood of $0\in S$.
\item An open set $\shU\subseteq \X$.
\item Constants $C_1$ and $C_2$
\end{itemize}
such that for $t$ sufficiently small,
$f^{-1}(t)\cap\shU$ is a deformation of $g^{-1}(t)$ of size
$C_1|t|^{C_2}$. Here $g^{-1}(t)$ is itself of the form $X(B,\alpha')$
for some affine structure on $B$ and \v Cech 1-cocycle $\alpha'$,
so it makes sense to talk about small deformations of $g^{-1}(t)$.
\end{theorem}

Thus, given an integral affine manifold with singularities $B$,
if we can find a toric polyhedral decomposition $\P$ of $B$,
we can construct $\X_0$, and hopefully a smoothing under nice
circumstances. Then $\X_t$ is a topological compactification of $X(U)$,
and we also have a good approximation to the complex structure
on an open subset of $\X_t$. Ideally, of course, we would like
to describe $\X_t$ explicitly as a topological compactification
of $X(U)$. However, this is not possible in general as there
may be a number of different compactifications. In special cases,
it is possible to prove stronger results, but we do not give details
here.

If $\X\rightarrow S$ is a polarized toric degeneration, with
degeneration data $(B,\P,\varphi)$, we obtain $(\check B,\check\P,\check
\varphi)$ via the discrete Legendre transform, and if we have a
corresponding polarized toric degeneration $\check\X\rightarrow S$,
then $\check\X_t$ is a compactification of some $X(\check B_0)$,
which is the dual torus bundle to $X(B_0)\rightarrow B_0$ by Proposition
\ref{intersectioncomplex}.
Hence we recover at least a topological form of SYZ, along with
additional information about complex structures. In addition, we
see mirror symmetry as a duality between affine manifolds via the 
discrete Legendre transform, clearly related to the continuous Legendre
transform appearing in SYZ.

Finally, let us underline the context of the error estimate in
terms of the ideas of \cite{GW}. If $t\in S^*$ with $r=-{\log |t|\over 
2\pi}\in (R,\infty)$, then $g^{-1}(t)\cong X(B,\alpha')$ where
$B$ has the affine structure given by $rd+d'$ and $\alpha'$
is some choice of twisting. One can alternatively rescale the affine
structure by multiplying by $\epsilon=1/r$, giving a developing map
$d+\epsilon d'$, which one should think of as a small perturbation
of $d$. Then $X(B,\alpha')$ can be
viewed as a twist of $\T_B/\epsilon\Lambda$, where $\Lambda\subseteq
\T_B$ is the lattice of integral vector fields coming from
the affine structure $d+\epsilon d'$. Thus we see as $t\rightarrow 0$,
$\epsilon\rightarrow 0$ and essentially the fibres are shrinking,
with radius $\epsilon$. The size of the deformation is then 
$O(e^{-C/\epsilon})$, i.e. decays exponentially in terms of the fibre
radius. This is a very similar picture to that of \cite{GW}. In fact,
if one has a potential function $\varphi$ on $U\subseteq B_0$ satisfying
the Monge-Amp\`ere equation, then one can obtain an almost Ricci-flat
metric on $f^{-1}(t)\cap\U$ for $t$ sufficiently small. As Ilia Zharkov has
also advocated \cite{Zh}, this is a first step towards proving a 
limiting form of the SYZ conjecture (see \cite{GExamples}). 
(In fact \cite{HZ2} contains a proof of the above theorem in the
toric hypersurface case).

\end{document}